\documentclass[11pt]{amsart}
\usepackage{amsfonts,amssymb,amsmath,kotex,mathtools,amsthm}

\usepackage{hyperref}

\usepackage{tikz-cd}
\usepackage{tikz}

\usepackage{enumerate} 
\usepackage{cellspace}
\usepackage{geometry}

\usepackage{stix}
\theoremstyle{definition}
\newtheorem{thm}{Theorem}[section]

\newtheorem{definition}[thm]{Definition}
\newtheorem{lemma}[thm]{Lemma}
\newtheorem{Cor}[thm]{Corollary}

\newtheorem{fact}[thm]{Fact}
\newtheorem{remark}[thm]{Remark}
\newtheorem{example}[thm]{Example}

\newtheorem*{ack}{Acknowledgements}
\newtheorem*{ai}{AI disclosures}
\theoremstyle{remark}

\title{Zariski dense hyperbolic $3$-manifold groups in $\mathrm{SL}(4,\mathbb{Z})$}

\author{DONG HYUN LEE}
\address{Department of Mathematical Sciences, Seoul National University, Seoul 08826,
South Korea}
\email{ldh0260@snu.ac.kr}

\begin{document}
\begin{abstract}
We show that $\mathrm{SL}(4,\mathbb{Z})$ contains an infinite family of Zariski-dense subgroups isomorphic to the fundamental group of a fixed closed hyperbolic $3$-manifold, which are pairwise noncommensurable up to conjugacy.
\end{abstract}
\maketitle

\section{Introduction}
Let $\Gamma$ be a finitely generated subgroup of $\mathrm{GL}(n,\mathbb{Z})$, and let $G=\mathrm{Zcl}(\Gamma)$ denote its Zariski closure in $\mathrm{GL}(n,\mathbb{R})$. The group $\Gamma$ is called a \emph{thin group} if it has infinite index in $G(\mathbb{Z})$ \cite{sarnak2014notes}. Despite their infinite index, thin groups often share important arithmetic and dynamical properties with lattices, such as strong approximation and expansion properties under suitable conditions \cite{salehigolsefidy2012expansion}. Nevertheless, constructing nonfree thin groups in $\mathrm{SL}(n,\mathbb{Z})$ remains difficult.

Higher Teichm\"uller theory has been used to produce thin groups; see \cite{10.2140/gt.2011.15.1,long2024zariski}. Long, Reid, and Thistlethwaite constructed thin groups isomorphic to surface groups by deforming representations within components of the representation variety consisting entirely of discrete and faithful representations. By combining these deformations with arithmetic arguments, they obtained examples of integral representations. In particular, they exhibited infinitely many discrete and faithful representations of the triangle group $\Delta(3,4,4)$ into $\mathrm{SL}(3,\mathbb{Z})$.

Our result continues in this line of research. We expand Long and Reid’s result to the existence of a fixed hyperbolic 3-manifold group in $\mathrm{SL}(4,\mathbb{Z})$. Examples of thin groups isomorphic to closed hyperbolic $3$-manifold groups in $\mathrm{SL}(4,\mathbb{Z})$ were already known; see, for example, \cite{ChoiChoi2015,audibert2025zariski}. However, to the author's knowledge, it was not previously known whether a fixed closed hyperbolic $3$-manifold group could admit infinitely many discrete and faithful representations into $\mathrm{SL}(4,\mathbb{Z})$ with Zariski-dense images.

\begin{thm}\label{1.1}
There exists a closed hyperbolic $3$-manifold $M$ such that $\pi_1(M)$ admits infinitely many discrete and faithful representations into $\mathrm{SL}(4,\mathbb{Z})$ with Zariski-dense images. Moreover, their images are pairwise noncommensurable up to conjugacy in $\mathrm{SL}(4,\mathbb{R})$.
\end{thm}

One possible approach is to deform the holonomy representation of a closed hyperbolic $3$-manifold in $\mathrm{SL}(4,\mathbb{R})$. Cooper, Long, and Thistlethwaite found experimentally that most of the $4,500$ two-generator closed hyperbolic $3$-manifolds they examined are infinitesimally projectively rigid \cite{cooper2006computing}. Thus, direct deformation of the holonomy representation provides few candidates. Moreover, even when nontrivial deformations exist, a positive-dimensional deformation space does not always provide infinitely many integral representations by itself.

Vinberg theory for Coxeter groups gives another approach. For an irreducible representation of a Coxeter group as a reflection group that is not virtually abelian, its Zariski closure is either $\mathrm{SL}^{\pm}(V)$ or an orthogonal group $\mathrm{O}_f(V)$ associated with a nondegenerate symmetric bilinear form $f$ \cite{audibert2025zariski}. Vinberg's theory therefore provides a flexible tool for constructing Zariski-dense representations. However, integrality remains a restriction. For compact hyperbolic Coxeter truncation polytopes of dimension at least $3$, Ko proved that the geometric component contains only finitely many integral representations up to conjugacy \cite{ko2026finiteness}. This shows that the main difficulty is not just to construct a deformation space, but to find infinitely many integral points in it.

Our construction combines Vinberg's theory with $\mathbb{H}^{p,q}$-convex cocompact deformation theory. We begin with an integral reflection representation of the Coxeter group associated with a compact hyperbolic Coxeter $3$-polytope, taking values in $\mathrm{O}(3,3)$. By the results of \cite{danciger2024convex}, this representation is $\mathbb{H}^{3,2}$-convex cocompact. We then use bending to produce infinitely many representations in the corresponding $\mathbb{H}^{3,2}$-convex cocompact component. The work of Beyrer and Kassel implies that these deformations remain discrete and faithful. Finally, we restrict the bending family to a fixed torsion-free subgroup of finite index, which is the fundamental group of a closed hyperbolic $3$-manifold.

The main arithmetic step is to choose infinitely many values of a bending parameter for which the resulting representations are integral. Using the exceptional isogeny between $\mathrm{SL}(4,\mathbb{R})$ and $\mathrm{SO}_0(3,3)$, together with the lifting and integrality arguments, we obtain infinitely many representations into $\mathrm{SL}(4,\mathbb{Z})$ with the same fundamental group. We can also verify that their images are Zariski-dense and represent infinitely many conjugacy classes.

In fact, the images of the representations in Theorem~\ref{1.1} are thin. The lattice $\mathrm{SL}(4,\mathbb{Z})$ has Kazhdan's property $(T)$, whereas the closed hyperbolic $3$-manifold group arising as a finite-index subgroup of an infinite Coxeter group does not \cite{bozejko1988infinite}. Since Kazhdan's property $(T)$ is preserved under passing to finite-index subgroups, the faithful images in Theorem~\ref{1.1} cannot have finite index in $\mathrm{SL}(4,\mathbb{Z})$.

It remains to distinguish their commensurability classes. Long and Reid proved a noncommensurability result for representations of a hyperbolic triangle group in the $\mathrm{SL}(3,\mathbb{R})$-Hitchin component \cite{10.1215/ijm/1498032023}. By the result of Choi and Goldman \cite{choi1993convex}, the $\mathrm{SL}(3,\mathbb{R})$-Hitchin component for a closed surface group coincides with the deformation space of convex real projective structures, and hence with the space of Benoist representations in this setting. On the other hand, Viaggi \cite{viaggi2025geometric} gave a unified proof of the results of Benoist and Koszul \cite{benoist2005convexes,koszul1968deformations} and of Beyrer and Kassel \cite{beyrer2025mathbb}, concerning higher higher Teichm\"uller components for Benoist representations and $\mathbb{H}^{p,q}$-convex cocompact representations, respectively. Therefore, this common viewpoint lead us to ask whether the noncommensurability phenomenon of Long and Reid extends to $\mathbb{H}^{p,q}$-convex cocompact representations.
\begin{thm}\label{1.2}
Let $\Gamma$ be the fundamental group of a closed hyperbolic $p$-manifold, where $p\geq 2$, and let $q\geq 1$. Suppose that $(\rho_n)_{n}$ is a sequence of Zariski-dense $\mathbb{H}^{p,q}$-convex cocompact representations $\rho_n:\Gamma\longrightarrow \mathrm{SO}(p,q+1)$ whose image subgroups are pairwise nonconjugate. Then there exists an infinite subsequence whose image subgroups are pairwise noncommensurable up to conjugacy.
\end{thm}

Applying Theorem~\ref{1.2} with $(p,q)=(3,2)$ to the infinite sequence in Theorem \ref{1.1} gives an infinite subsquence whose images are pairwise noncommensurable up to conjugacy, and hence completes the proof of Theorem~\ref{1.1}.

\begin{ai}
The author used ChatGPT-6 for drawing figures and for proofreading the introduction at the end of the writing process. All mathematical content was produced without the aid of AI.
\end{ai}

\begin{ack}
I am really grateful to my advisor, Gye-Seon Lee for suggesting this problem and sharing his idea. I would like to thank Jacques Audibert for his valuable insights, which I learned about through my advisor. I also thank Seong Hwan Ko and Seong Youn Kim for helpful comments on the readability of proofs. This work was supported by the National Research Foundation of Korea(NRF) grant funded by the Korea government(MSIT) (RS-2025-02293115).
\end{ack}
\section{Preliminaries}
\subsection{\texorpdfstring{$\mathbb{H}^{p,q}$}{H(p,q)} convex-cocompact representations}
Our exposition in this subsection follows closely \cite{beyrer2025mathbb,viaggi2025geometric}. A subset of projective space $\mathbb{P}(\mathbb{R}^n)$ is said to be \textit{convex} if there exists an affine chart containing the subset such that it is convex in the corresponding affine space. In particular, a convex subset in $\mathbb{P}(\mathbb{R}^n)$ is called \textit{properly convex} if it does not contain any affine line.

\begin{definition}
    A discrete subgroup of $\Gamma$ of $\mathrm{PGL}(V)$ is said to be \textit{(strongly) convex cocompact} in $\mathbb{P}(V)$ if the following conditions hold:
    \begin{enumerate}
        \item it acts properly discontiuously on some properly convex open subset of $\mathbb{P}(V)$,
        \item it acts cocompactly on some nonempty closed convex subset $\mathcal{C}$ of $\Omega$,
        \item $\Omega$ can be taken so that its boundary $\partial\Omega:=\overline{\Omega}\setminus\Omega$ does not contain any nontrivial projective segment and every point has a unique supporting hyperplane.
    \end{enumerate}
\end{definition}
Here $\overline{\Omega}$ denotes the closure of $\Omega$ in $\mathbb{P}(V)$.

\hspace{1em}

This notion can be extended to a pseudo-Riemannian hyperbolic space by \cite{danciger2018convex}. We equip $\mathbb{R}^{p+q+1}$ with the symmetric bilinear form
\[\langle x,y\rangle_{p,q+1}:=x_1y_1+\ldots + x_py_p-x_{p+1}y_{p+1}-\ldots -x_{p+q+1}y_{p+q+1},\]
and denote this space by $\mathbb{R}^{p,q+1}$.

\begin{definition} The \textit{pseudo-Riemannian
hyperbolic space} $\mathbb{H}^{p,q}$ is
\[\mathbb{H}^{p,q}=\{x\in \mathbb{R}^{p,q+1}:\langle x,x\rangle_{p,q+1}=-1\}\]
equipped with the pseudo-Riemannian metric $g_{\mathbb{H}^{p,q}}$ of signature $(p,q)$ induced by
the restriction of the symmetric bilinear form $\langle \cdot,\cdot\rangle_{p,q+1}$ to the tangent spaces $T_x\mathbb{H}^{p,q} = x^{\perp}$.
\end{definition}
In projective geometry, we can consider the pseudo-Riemannian hyperbolic space $\mathbb{H}^{p,q}$ as a space of the projective space $\mathbb{P}(\mathbb{R}^{p,q+1})$
\[\mathbb{H}^{p,q}=\{[x]\in \mathbb{P}(\mathbb{R}^{p,q+1}):\langle x,x\rangle_{p,q+1}<0\}.\]

\begin{definition}
A discrete subgroup $\Gamma$ of $\mathrm{PO}(p,q+1)$ is $\mathbb{H}^{p,q}$\textit{-convex cocompact} if it acts properly discontinuously and cocompactly on some properly convex closed subset $\mathcal{C}$ of $\mathbb{H}^{p,q}$ with nonempty interior whose ideal
boundary $\partial_i\mathcal{C}:=\overline{\mathcal{C}}\setminus\mathcal{C}$ does not contain any nontrivial projective segment.
\end{definition}
Here, we consider $\mathcal{C}$ as the closure of $\mathcal{C}$ in $\mathbb{P}(\mathbb{R}^{p,q+1})$.

\begin{definition}
    The \textit{cohomological dimesion} of a torsion-free group $\Gamma$ is defined as the largest number $n$ in $\mathbb{N}$ which satisfies $H^n(\Gamma)\neq0$ for some $\mathbb{Z}\Gamma$-module $M$, and is denoted by $\mathrm{cd}(\Gamma)$.
\end{definition}
By Serre’s theorem, if $
\Gamma$ is a torsion-free group and $\Gamma'$ is a finite-index subgroup of $\Gamma$, then $\mathrm{cd}(\Gamma) = \mathrm{cd}(
\Gamma')$. Consequently, if a group $\Gamma$ admits torsion-free finite-index subgroups $\Gamma_1$ and $\Gamma_2$, their cohomological dimensions coincide; specifically, $\mathrm{cd}(\Gamma_1) = \mathrm{cd}(\Gamma_1\cap \Gamma_2) = \mathrm{cd}(\Gamma_2)$. This gives the well-defined notion of the \textit{virtual cohomological dimension of $\Gamma$}, denoted by $\mathrm{vcd}(\Gamma)$, for any group having a torsion-free finite-index subgroup.

The following fact implies that the set $\mathrm{Hom}^{\text{cc}}_{0}(\Gamma,\mathrm{O}(p,q+1))$  of discrete and faithful $\mathbb{H}^{p,q}$-convex cocompact representations of $\Gamma$ is a union of connected components of $\mathrm{Hom}(\Gamma,\mathrm{O}(p,q+1))$.

\begin{fact}[{\cite[Thm 7.7]{beyrer2025mathbb}}]\label{thm:HpqCC}
    Let $p \geq 2$ and $q \geq 1$ be integers, and let $\Gamma$ be a Gromov hyperbolic group with
    $\mathrm{vcd}(\Gamma) = p$. Then the set of $\mathbb{H}^{p,q}$-convex cocompact representations of $\Gamma$ is a union of connected
components of $\mathrm{Hom}(\Gamma,\mathrm{PO}(p,q +1))$. In particular, for any finite normal subgroup $\Gamma'$ of $\Gamma$, the set $\mathrm{Hom}_{\Gamma'}^{\text{cc}}(\Gamma,\mathrm{PO}(p,q+1))$ of $\mathbb{H}^{p,q}$-convex cocompact representations with kernel exactly $\Gamma'$ is a union of connected components of $\mathrm{Hom}(\Gamma,\mathrm{PO}(p,q+1))$.
\end{fact}
To prove Theorem 1.2, we need the notion of \textit{complete maximal $p$-submanifold} in $\mathbb{H}^{p,q}$.
\begin{definition}
    A smooth $p$-submanifold $M$ in $\mathbb{H}^{p,q}$ is \textit{spacelike} if the restriction of $g_{\mathbb{H}^{p,q}}$ to $T_xM$ is positive definite for every $x\in M$. In particular, $g_{\mathbb{H}^{p,q}}$ induces a Riemannian metric on $M$. The spacelike $p$-submanifold is said to be \textit{complete} if it is complete in the Riemannian sense.
\end{definition}

\begin{definition}A spacelike $p$-submanifold $M$ in $\mathbb{H}^{p,q}$ is \textit{maximal} if the second fundamental form $\mathrm{II}_M$ of $M$ satisfies the following equation $\mathrm{tr}_M(\mathrm{II}_M)=0$.
\end{definition}

\begin{fact}[{\cite[Cor 1.1]{seppi2023complete}}]\label{complete_maximal_submanifold}
     Let $\Gamma$ be a Gromov group with Gromov boundary homeomorphic to a $(p-1)$-sphere. Every $\mathbb{H}^{p,q}$-convex-cocompact representation $\rho:\Gamma\rightarrow \mathrm{PO}(p,q+1)$ preserves a unique complete
maximal $p$-submanifold $M$ in $\mathbb{H}^{p,q}$. Furthermore
\begin{enumerate}[(1)]
    \item $\rho$ acts properly discontinuously and cocompactly on $M$,
    \item $M$ depends real analytically on $\rho$.
\end{enumerate}

In fact, if $\Gamma$ is an infinite discrete subgroup of $\mathrm{PO}(p,q+1)$ that is $\mathbb{H}^{p,q}$-convex cocompact, then $\Gamma$ is Gromov hyperbolic and has virtual cohomological dimension at most $p$. Moreover, $\mathrm{vcd}(\Gamma)=p$ if and only if the Gromov boundary of $\Gamma$ is homeomorphic to a $(p-1)$-dimensional sphere; see \cite[Cor 11.10]{DGK24}.

\end{fact}
\subsection{Vinberg theory}
In this subsection, we introduce some notions from Vinberg theory that will be used in the construction of thin groups. Our discussion follows closely \cite{audibert2025zariski,danciger2024convex}.
\begin{definition} A \textit{Coxeter matrix} $M$ on a finite set $S$ is a symmetric $|S|\times |S|$ matrix
$M=(m_{st})_{s,t\in S}$ with entries $m_{st}\in \{1,2,\ldots ,\infty \}$ with diagonal entries $m_{ss} =1$ and off-diagonal
entries $m_{st}\neq 1$. To a Coxeter matrix $M$ is associated a Coxeter group $W_S$: the group presented by the set of generators $S$ and the relations $(st)^{m_{st}} = 1$ for each $(s,t)\in S\times S$ with $m_{st}\neq \infty$.
\end{definition}
A (abstract) Coxeter group $W_S$ associated with a Coxeter matrix $M$ is defined by the presentation
\[W_S =\langle S |(st)^{m_{st}} = 1 \text{ for all } s,t\in S \text{ with } m_{st} \neq \infty \rangle.\]
The rank of $W_S$ is defined as the cardinality $|S|$ of $S$.
\begin{definition}
A \textit{Coxeter diagram} $\mathscr{G}_{W_S}$ of $W_S$ is an edge-labeled graph. The vertex set of $\mathscr{G}_{W_S}$ is $S$. Two vertices $s, t \in S$ are joined by an edge $\overline{st}$ if and only if $m_{st} \in \{3,4, \ldots ,\infty \}$, and the label of the edge $\overline{st}$ is $m_{st}$.
\end{definition}

\begin{definition}
    A Coxeter group $W_S$ is said to be \textit{irreducible} if its Coxeter diagram $\mathscr{G}_{W_S}$ is connected.
\end{definition}
Margulis and Vinberg classified irreducible Coxeter groups into three classes.
\begin{definition}[{\cite[Thm 1]{MV00}}] Every irreducible Coxeter group $W_S$ belongs to one of the followings:
\begin{enumerate}
    \item \textit{spherical} if $W_S$ is finite,
    \item \textit{affine} if $W_S$ is virtually abelian and infinite,
    \item \textit{large} if $W_S$ has a finite index subgroup with a non-abelian free quotient. 
\end{enumerate}
\end{definition}
\begin{remark}
    The irreducible spherical and irreducible affine Coxeter groups were classified by Coxeter \cite{Cox32,Cox34}.
\end{remark}
We produce discrete and faithful reprsentations from data of a Coxeter diagram.
\begin{definition}
A \textit{Cartan matrix} $\mathscr{A}_S=(\mathscr{A}_{ij})_{i,j\in S}$ on $S$ is a $|S|\times |S|$ matrix satisfying the following conditions:
\begin{enumerate}
    \item $\mathscr{A}_{ii}=2$ for all $i\in S$,
    \item $\mathscr{A}_{ij}=0$ for all $i\neq j$ with $m_{ij}=2$,
    \item $\mathscr{A}_{ij}<0$ for all $i\neq j$ with $m_{ij}\neq 2$,
    \item $\mathscr{A}_{ij}\mathscr{A}_{ji}=4\cos^2(\pi/m_{ij})$ for all $i\neq j$ with $m_{ij}\in\{3,4,\ldots\}$.
\end{enumerate}
Such a Cartan matrix is said to be \textit{weakly compatible} with $W_S$. It is said to be \textit{compatible} with $W_S$ if, in addition,
$\mathscr{A}_{ij}\mathscr{A}_{ji}\geq 4$ for all $i\neq j$ with $m_{ij}=\infty$.
\end{definition}

\begin{definition}
    Let $W_S$ be a Coxeter group. A representation $\rho:W_S \rightarrow \mathrm{GL}(V)$ is \textit{a representation of $W_S$ as a reflection group} if for every $s\in S$ there exist $v_s\in V$ and $\alpha_s \in V^*$ such that:
    \begin{enumerate}
        \item for every $s\in S$, $\rho(s)=\mathrm{Id}-\alpha_s\otimes v_s$,
        \item $\mathscr{A} =(\alpha_s(v_t))_{s,t\in S}$ is a Cartan matrix associated with $W_S$,
        \item  the Cartan matrix $\mathscr{A}$ and the Coxeter group $W_S$ are compatible,
        \item the convex cone $\widetilde{\Delta}:=\bigcap_{s\in S}\{v\in V:\alpha_s(v)\leq 0\}$ has nonempty interior.
    \end{enumerate}
 We denote by $\mathrm{Hom}^{\text{ref}}(W_S,\mathrm{GL}(V))$ the set of all representations of $W_S$ as a reflection group in $V$. 
\end{definition}

The pair $(\alpha_s,v_s)$ defining a reflection is unique up to the action 
\[\lambda\cdot(\alpha_s,v_s):=(\lambda\alpha_s,\lambda^{-1}v_s)\]
for $\lambda \in \mathbb{R}^{\times}$. Thus, the reflections $\rho(s_1),\ldots, \rho(s_{|S|})$ determine the $|S|\times |S|$ Cartan matrix  up to conjugation by a diagonal matrix. Condition $(3)$ in the definition of a Cartan matrix guarantees that the diagonal matrix can be chosen to be positive.
\begin{definition}
Two $|S|\times |S|$ Cartan matrices associated with $W_S$ are considered \textit{equivalent} if they differ by conjugation by a positive diagonal matrix.
\end{definition}
In particular, we shall say that $\rho:W_S\rightarrow \mathrm{GL}(V)$ is \textit{symmetrizable} if the Cartan matrix $\mathscr{A}$ is equivalent to a symmetric matrix.
\begin{example}\label{geom_repn}
Let $W_S$ be a Coxeter group. To any Cartan matrix $\mathscr{A}$ compatible with $W_S$, one can associate a representation $\rho_{\mathscr{A}}:W_S\rightarrow \mathrm{GL}(|S|,\mathbb{R})$ defined by sending each $s\in S$ to the reflection
\[\rho_{\mathscr{A}}(s):=(v\mapsto v-(v^T\mathscr{A}e_s)e_s).\]
\end{example}

When $W_S$ is irreducible and large, these reflection representations form a union of connected components in the space of discrete and faithful representations. More precisely, we have the following fact.

\begin{fact}
Suppose that $W_S$ is irreducible and large. Then $\mathrm{Hom}^{\mathrm{ref}}(W_S,\mathrm{GL}(n,\mathbb{R}))$ is both open and closed in the space $\mathrm{Hom}^{\mathrm{fd}}(W_S,\mathrm{GL}(n,\mathbb{R}))$ of discrete and faithful representations of $W_S$ into $\mathrm{GL}(n,\mathbb{R})$.
\end{fact}

\begin{fact}
     Let $W_S$ be a Coxeter group and $\rho:W_S\rightarrow \mathrm{GL}(V)$ a representation of $W_S$ as a reflection group. Then:
\begin{enumerate}
    \item the representation $\rho$ is faithful and $\rho(W_S)$ is discrete in $\mathrm{GL}(V)$;
    \item the union of the convex cones $\rho(\gamma)(\widetilde{\Delta})$ as $\gamma$ varies within $W_S$ is a convex cone $\mathscr{C}_{TV}$ of $V$;
    \item the action of $\rho(W_S)$ on the interior $\Omega_{TV}$ of $\mathbb{P}(\mathscr{C}_{TV})$ is proper. The domain $\Omega_{TV}$ is called the Tits–Vinberg domain of the reflection group $\rho(W_S)$.
\end{enumerate} 
\end{fact}

The following dichotomy provides a convenient criterion for determining whether the image of a reflection representation of a Coxeter group is Zariski-dense.
\begin{fact}[\cite{audibert2025zariski}]\label{fact:zariski_dense}
    Let $W$ be a finitely generated Coxeter group that is not virtually abelian and $\rho:W \rightarrow \mathrm{GL}(V)$ a representation of $W$ as a reflection group. Suppose that $\rho$ is irreducible.
    \begin{enumerate}[(1)]
        \item If $\rho$ preserves a nonzero symmetric bilinear form $f$ on $V$, then the Zariski-closure of $\rho(W)$ is the orthogonal group $\mathrm{O}_f(V)$ of $f$.
        \item Otherwise, the Zariski-closure of $\rho(W)$ is $\mathrm{SL}^{\pm}(V)$.
    \end{enumerate}
\end{fact}
Irreducibility of a reflection representation of an irreducible Coxeter group can be checked by using the criterion below. 
\begin{fact}\label{irred} Let $W_S$ be an irreducible Coxeter group and $\rho:W_S \rightarrow \mathrm{GL}(V)$ a representation of $W_S$ as a reflection group. Let $\mathscr{A}$ be the Cartan matrix of $\rho$. Then:
\begin{enumerate}[(1)]
    \item The representation $\rho$ is irreducible if and only if $\mathrm{dim}(V)=\mathrm{rank}(\mathscr{A})$.
    \item Assume $\rho$ is irreducible. Then $\rho$ preserves a nonzero (hence nondegenerate) symmetric bilinear form on $V$ if and only if $\mathscr{A}$ is symmetrizable.
\end{enumerate} 
\end{fact}

\begin{example}\label{example1}Let $W_{23}$ be the Coxeter group of cu23 in  \cite{choi2012projective} with $6$ generators for which the Coxeter diagram is given by Figure \ref{fig:cube-and-coxeter-diagram}.
This Coxeter group is in fact large from the classification of irreducible spherical and irreducible affine Coxeter groups.
\begin{figure}[htbp]
\centering
\begin{tikzpicture}

\begin{scope}[
  line width=0.7pt,
  edge label/.style={font=\normalsize,inner sep=2pt},
  face/.style={font=\large}
]
  \coordinate (A) at (-2, 2);
  \coordinate (B) at ( 2, 2);
  \coordinate (C) at ( 2,-2);
  \coordinate (D) at (-2,-2);

  \coordinate (a) at (-0.65, 0.65);
  \coordinate (b) at ( 0.65, 0.65);
  \coordinate (c) at ( 0.65,-0.65);
  \coordinate (d) at (-0.65,-0.65);

  \draw (A) -- node[edge label,above] {$2$} (B);
  \draw (B) -- node[edge label,right] {$3$} (C);
  \draw (C) -- node[edge label,below] {$2$} (D);
  \draw (D) -- node[edge label,left]  {$2$} (A);

  \draw (a) -- node[edge label,above] {$3$} (b);
  \draw (b) -- node[edge label,right] {$2$} (c);
  \draw (c) -- node[edge label,below] {$3$} (d);
  \draw (d) -- node[edge label,left]  {$2$} (a);

  \draw (A) -- node[edge label,above right] {$3$} (a);
  \draw (B) -- node[edge label,above left]  {$2$} (b);
  \draw (C) -- node[edge label,below left]  {$2$} (c);
  \draw (D) -- node[edge label,below right] {$3$} (d);

  \node[face] at ( 0.00, 0.00) {$F_1$};
  \node[face] at (-1.32, 0.00) {$F_2$};
  \node[face] at ( 1.32, 0.00) {$F_3$};
  \node[face] at ( 2.75, 0.00) {$F_4$};
  \node[face] at ( 0.00, 1.32) {$F_5$};
  \node[face] at ( 0.00,-1.32) {$F_6$};
\end{scope}

\begin{scope}[
  xshift=6.8cm,
  yshift=-0.6cm,
  line width=0.8pt,
  vertex/.style={
    circle, draw, fill=white,
    inner sep=0pt, minimum size=8mm, font=\large
  },
  edge label/.style={font=\large,inner sep=2pt}
]
  \node[vertex] (s1) at (-2.1, 0.0) {$s_1$};
  \node[vertex] (s2) at ( 2.1, 0.0) {$s_2$};
  \node[vertex] (s3) at ( 2.1, 2.2) {$s_3$};
  \node[vertex] (s4) at (-2.1, 2.2) {$s_4$};
  \node[vertex] (s5) at ( 0.0, 1.0) {$s_5$};
  \node[vertex] (s6) at ( 0.0,-1.0) {$s_6$};

  \draw (s1) -- node[edge label,above left]  {$3$} (s5);
  \draw (s1) -- node[edge label,below left]  {$3$} (s6);
  \draw (s2) -- node[edge label,above right] {$3$} (s5);
  \draw (s2) -- node[edge label,below right] {$3$} (s6);
  \draw (s3) -- node[edge label,above]       {$3$} (s4);

  \draw (s1) -- node[edge label,left]  {$\infty$} (s4);
  \draw (s2) -- node[edge label,right] {$\infty$} (s3);
  \draw (s5) -- node[edge label,right] {$\infty$} (s6);
\end{scope}

\end{tikzpicture}
\caption{A Coxeter cube cu23 and its Coxeter diagram.}
\label{fig:cube-and-coxeter-diagram}
\end{figure}
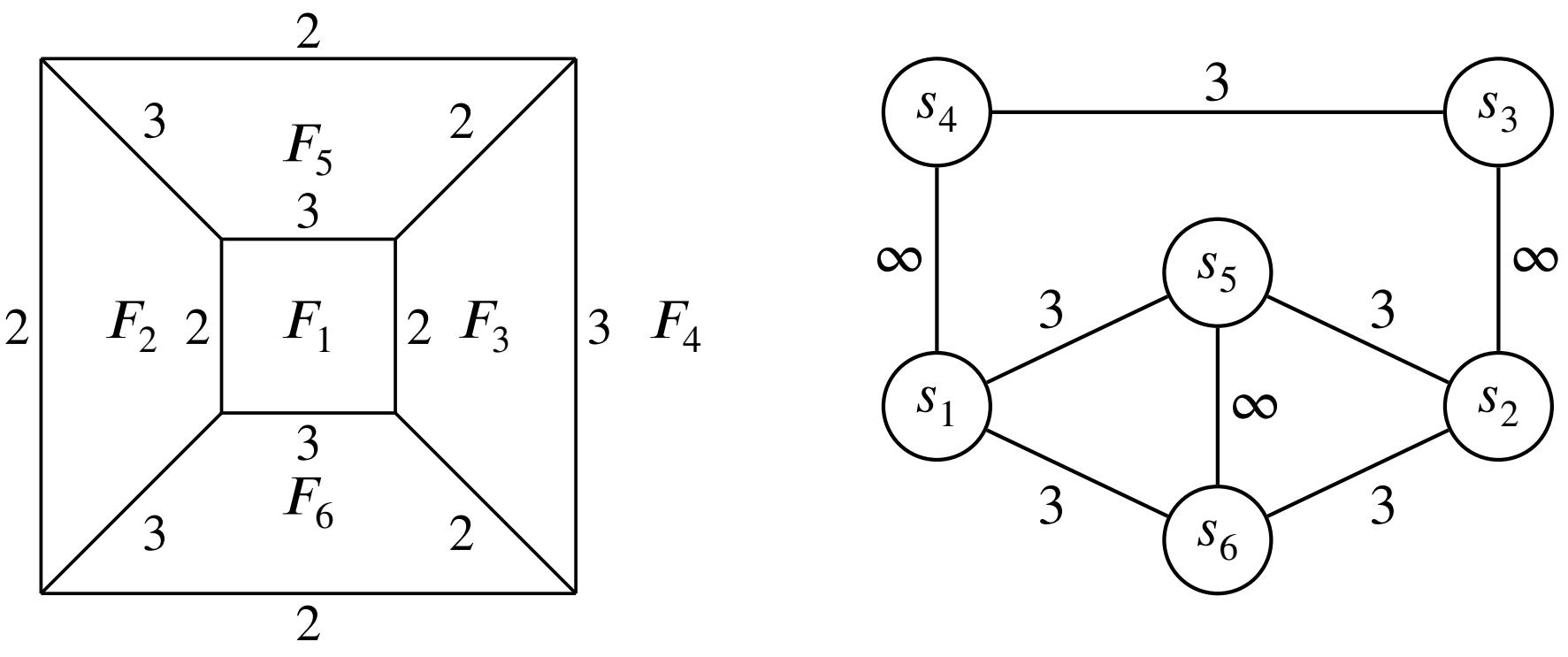
Its symmetrizable Cartan matrix $\mathscr{A}$ for $W_S$ can be conjugated by diagonal matrices to 
\begin{align*}
        \mathscr{A} = \begin{bmatrix*}[r]
            2 & 0 & 0 & -u & -1 & -1 \\
            0 & 2 & -v & 0 & -1 & -1 \\
            0 & -v & 2 & -1 & 0 & 0 \\
            -u & 0 & -1 & 2 & 0 & 0 \\
            -1 & -1 & 0 & 0 & 2 & -w \\
            -1 & -1 & 0 & 0 & -w & 2   
        \end{bmatrix*}
\end{align*}
for some $u,v,w\geq 2$. The Fact $\ref{irred}$ implies that if $\mathrm{det}(\mathscr{A})\neq 0$ then the representation $\rho$ induced by $\mathscr{A}$ is irreducible. Since its $\{1,2,5,6\}$-principal minor matrix  
\begin{align*}
        \mathscr{A}_{1256} := \begin{bmatrix*}[r]
            2 & 0 & -1 & -1 \\
            0 & 2& -1 & -1 \\
            -1 & -1 & 2 & -w \\
            -1 & -1 & -w & 2   
        \end{bmatrix*}
\end{align*}
is of signature $(3,1)$ for every $w\geq 2$, a linear algebra calculation shows that if $\mathrm{det}(\mathscr{A})<0$ then $\mathscr{A}$ is of signature $(5,1)$ or $(3,3)$. In this case, its Zariski closure of $\rho(W_S)$ associated with $\mathscr{A}$ is conjugate to either $\mathrm{O}(5,1)$ or $\mathrm{O}(3,3)$.
\end{example}

\begin{Cor}
Suppose that $W_S$ is a large, word-hyperbolic, irreducible Coxeter group, and $n=p+q+1\geq |S|$. Let $\mathrm{Hom}_0^{\text{cc}}(W_S,\mathrm{O}(p,q+1))$ be the set of discrete and faithful $\mathbb{H}^{p,q}$-convex cocompact representations. Every component containing a representation of $W_S$ as a reflection group in $\mathrm{Hom}^{\text{cc}}_0(W_S,\mathrm{O}(p,q+1))$ is contained in $\mathrm{Hom}^{\text{ref}}(W_S,\mathrm{GL}(n,\mathbb{R}))$ where $n=p+q+1$. 
\end{Cor}
\begin{proof}
    Let $\mathscr{C}$ be a connected component of $\mathrm{Hom}^{\text{cc}}_{0}(W_S,\mathrm{O}(p,q+1))$ containing a representation of $W_S$ as a reflection group. Note that $i_*(\mathscr{C})$ is connected in $\mathrm{Hom}^{\text{fd}}(W_S,\mathrm{GL}(n,\mathbb{R}))$ where $i_*:\mathrm{Hom}(W_S,\mathrm{O}(p,q+1))\rightarrow \mathrm{Hom}(W_S,\mathrm{GL}(n,\mathbb{R}))$ is induced by a natural embedding $i:\mathrm{O}(p,q+1)\rightarrow \mathrm{GL}(n,\mathbb{R})$. Since $\mathrm{Hom}^{\text{ref}}(W_S,\mathrm{GL}(n,\mathbb{R}))$ is open and closed in $\mathrm{Hom}^{\text{fd}}(W_S,\mathrm{GL}(n,\mathbb{R}))$, we see that $i_*(\mathscr{C})\subseteq \mathrm{Hom}^{\text{ref}}(W_S,\mathrm{GL}(n,\mathbb{R})).$
\end{proof}

We denote by $\mathbf{X}(\Gamma, G)$ the (real GIT) quotient of the space of representations from $\Gamma$ to a real Lie group $G\leq \mathrm{GL}(n, \mathbb{R})$ under the conjugation action of $\mathrm{GL}(n, \mathbb{R})$, defined as
\begin{align*}
\mathbf{X}(\Gamma, G) :&= (\pi \circ i_*)(\mathrm{Hom}(\Gamma, G))
\end{align*}
where $i_*:\mathrm{Hom}(\Gamma,G)\rightarrow \mathrm{Hom}(\Gamma,\mathrm{GL}(n,\mathbb{R}))$ is induced by a natural embedding $i:G\rightarrow \mathrm{GL}(n,\mathbb{R})$ and the quotient map $\pi:\mathrm{Hom}(\Gamma,\mathrm{GL}(n,\mathbb{R}))\rightarrow \mathrm{Hom}(\Gamma,\mathrm{GL}(n,\mathbb{R}))/\!/\mathrm{GL}(n,\mathbb{R})$, see \cite{RS90,BL21} for the real GIT quotient. Richardson and Slodowy proved that the closure of any orbit contains exactly a closed orbit corresponding to a semisimple representation \cite{RS90}. Thus, we can also consider $\mathbf{X}(\Gamma,G)$ as follows:
\[\mathbf{X}(\Gamma, G) = (\pi\circ i_*)(\mathrm{Hom}^{\text{ss}}(\Gamma, G)).\]
If $\Gamma$ is a Coxeter group $W_S$, then the quotient of semisimple representations of $W_S$ as a reflection group in $\mathbb{R}^n$ is denoted by 
\[\mathbf{X}^{\text{ref}}(W_S,\mathrm{GL}(n,\mathbb{R})):=(\pi\circ i_*)(\mathrm{Hom}^{\text{ref}}(W_S,\mathrm{GL}(n,\mathbb{R})).\]

Now we can parametrize the quotient space of semisimple $\mathbb{H}^{p,q}$-convex cocompact reflection representations with Cartan matrices of $W_S$ via the following facts.
\begin{fact}[\cite{danciger2024convex}]
The map assigning to a conjugacy class of semisimple representations the equivalence class of its Cartan matrices is a homeomorphism from the open and closed subset of $\mathbf{X}(W_S,\mathrm{GL}(V))$ consisting of representations generated by weakly compatible reflections to the space of equivalence classes of $|S|\times |S|$ matrices of rank $\leq \mathrm{dim}(V)$ that are weakly compatible with $W_S$. In particular, $\mathbf{X}^{\text{ref}}(W_S,\mathrm{GL}(V))$ is homeomorphically parametrized by the space of equivalence classes of $|S|\times |S|$ Cartan matrices of rank $\leq \mathrm{dim}(V)$ that are compatible with $W_S$. 
\end{fact}

\begin{fact}[\cite{danciger2024convex}]
Let $W_S$ be an infinite word hyperbolic irreducible Coxeter group. For any $V$ with $\mathrm{dim} V \geq |S|$, the set of characters $[\rho]\in \mathbf{X}^{\text{ref}}(W_S,\mathrm{GL}(V))$ for which $\rho(W_S)$ is convex cocompact in $\mathbb{P}(V)$ is precisely the interior of $\mathbf{X}^{\text{ref}}(W_S,\mathrm{GL}(V))$ in $\mathbf{X}(W_S,\mathrm{GL}(V))$. Moreover, this interior is homeomorphically parametrized by the equivalence class of its $|S|\times |S|$ Cartan matrices of rank$\leq \mathrm{dim} V$ that are weakly compatible with $W_S$ and $\mathscr{A}_{ij}\mathscr{A}_{ji}>4$ for $i\neq j$ with $m_{ij}=\infty$.
\end{fact}

\begin{thm}[Folklore]\label{thm:folklore}
Let $W_S$ be a large, word-hyperbolic, irreducible Coxeter group, and let $n=p+q+1\geq |S|$. Then the quotient space of semisimple $\mathbb{H}^{p,q}$-convex cocompact reflection representations of $W_S$ on $\mathbb{R}^n$ is homeomorphic to the space of equivalence classes of symmetrizable Cartan matrices $\mathscr{A}$ that are compatible with $W_S$, have signature $(p',q'+1)$ for some $p'\leq p$ and $q'\leq q$, and satisfy $\mathscr{A}_{ij}\mathscr{A}_{ji}>4$ for all $i\neq j$ with $m_{ij}=\infty$.
\end{thm}
\begin{proof}
Let $\mathscr{C}$ be a connected component of $\mathrm{Hom}^{\text{cc}}_{0}(W_S,\mathrm{O}(p,q+1))$ containing a representation of $W_S$ as a reflection group. Note that $i_*(\mathscr{C})$ is connected in $\mathrm{Hom}^{\text{fd}}(W_S,\mathrm{GL}(n,\mathbb{R}))$. 
Since $\mathrm{Hom}^{\text{ref}}(W_S,\mathrm{GL}(n,\mathbb{R}))$ is open and closed in $\mathrm{Hom}^{\text{fd}}(W_S,\mathrm{GL}(n,\mathbb{R}))$, we see that 
\[i_*(\mathscr{C})\subseteq \mathrm{Hom}^{\text{ref}}(W_S,\mathrm{GL}(n,\mathbb{R})).\] 
Thus, it follows that
\[\pi(i_*(\mathrm{Hom}^{\text{cc}}_{0}(W_S,\mathrm{O}(p,q+1))\cap \mathrm{Hom}^{\text{ref}}(W_S,\mathrm{GL}(n,\mathbb{R})))\subseteq \mathrm{int}(\mathbf{X}^{\text{ref}}(W_S,\mathrm{GL}(n,\mathbb{R}))).\]
Hence the previous facts conclude this proof. 
\end{proof}

\subsection{Algebraic bending deformation}
Let $\mathcal{O}$ be an $n$-dimensional orbifold, and let $\Sigma$ be an embedded, two-sided connected hypersurface that separates $\mathcal{O}$ into two pieces $\mathcal{O}_1$ and $\mathcal{O}_2$. Then $\pi_1^{\text{orb}}(\mathcal{O})$ is the amalgamated product $ \pi_1^{\text{orb}}(\mathcal{O}_1)*_{\pi_1^{\text{orb}}(\Sigma)}\pi_1^{\text{orb}}(\mathcal{O}_2)$. Let $\rho\in \mathrm{Hom}(\pi_1^{\text{orb}}(\mathcal{O}),G)$ be a representation. Suppose that $g\in  C_G(\rho(\pi_1^{\text{orb}}(\Sigma)))$ but does not commute with $\rho(\pi_1^{\text{orb}}(\mathcal{O}_2))$. Then we obtain another representation $\rho_g\in \mathrm{Hom}(\pi_1^{\text{orb}}(\mathcal{O}),G)$ by defining
\begin{align*}
    \rho_g(\gamma)=
    \begin{dcases*}
        \rho(\gamma) & for $\gamma\in \pi_1^{\mathrm{orb}}(\mathcal{O}_1)$\\
        g\rho(\gamma)g^{-1} & for $\gamma\in \pi_1^{\mathrm{orb}}(\mathcal{O}_2)$.
    \end{dcases*}
\end{align*}

\section{construction}
In this section we prove the main theorem \ref{1.1}. It is well known that there is an exceptional $\mathbb{Q}$-isogeny between $\mathrm{SL}(4,\mathbb{R})$ and $\mathrm{SO}(3,3)_0$, and thus $\mathrm{SL}(4,\mathbb{Z})$ is commensurable to $\mathrm{O}_{\mathbb{Z}}(3,3)$ by the following fact.

\begin{fact}[\cite{morris2015introductionarithmeticgroups}]
Suppose that $\varphi: G \rightarrow G'$ is a surjective homomorphism defined over $\mathbb{Q}$, and $G_{\mathbb{Z}}$ denotes the intersection $G \cap \mathrm{GL}(n, \mathbb{Z})$ (resp., $G'_{\mathbb{Z}}$). Then $\varphi(G_{\mathbb{Z}})$ is commensurable with $G'_{\mathbb{Z}}$.
\end{fact}
Via conjugation by $P \in \mathrm{GL}(6, \mathbb{Q})$, the orthogonal group $\mathrm{O}_f(\mathbb{Z})$ preserving
\[P \mathrm{diag}(1, 1, 1, -1, -1, -1)P^T\]
is commensurable with $\mathrm{O}_{\mathbb{Z}}(3, 3)$. In this context, we introduce an equivalence relation between two symmetric matrices.

\begin{definition}
    Two symmetric matrices(equivalently, quadratic forms) $Q_1$ and $Q_2$ over a field $k$ are \textit{$k$-equivalent} if there is $P\in \mathrm{GL}(n,k)$ such that $Q_2=PQ_1P^{T}$, and we denote it by $Q_1\simeq_{k}Q_2$.
\end{definition}
Now we shall construct a representation  $\rho_{\mathscr{A}}$ in $\mathrm{Hom}^{\text{fd}}(W_S,\mathrm{O}_f(\mathbb{Z}))$ by using Vinberg's construction for a Coxeter group $W_S$. From this representation, an infinite family of $(\rho_n)$ can be obtained via bending deformation. For the concrete computations, see Appendix.
\begin{lemma}
    There exists a closed hyperbolic 3-orbifold $\mathcal{O}$ such that its orbifold fundamental group $\pi_1^{\mathrm{orb}}(\mathcal{O})$ admits infinitely many discrete faithful representations into $\mathrm{O}_f(\mathbb{Z})$ where $\mathrm{O}_f(\mathbb{R})$ is $\mathbb{Q}$-equivalent to $\mathrm{O}(3,3)$.
\end{lemma}
\begin{proof}Let $W_{23}$ and $\mathscr{A}$ be the Coxeter group of hyperbolic 3-orbifold  cu23 and its Cartan matrix in example \ref{example1}. Substituting $(u, v, w) = (4, 3, 34)$ yields the matrix
    \begin{align*}
        \mathscr{A} = \begin{bmatrix*}[r]
            2 & 0 & 0 & -4 & -1 & -1 \\
            0 & 2 & -3 & 0 & -1 & -1 \\
            0 & -3 & 2 & -1 & 0 & 0 \\
            -4 & 0 & -1 & 2 & 0 & 0 \\
            -1 & -1 & 0 & 0 & 2 & -34 \\
            -1 & -1 & 0 & 0 & -34 & 2   
        \end{bmatrix*}
    \end{align*}
    This Cartan matrix $\mathscr{A}$ is $\mathbb{Q}$-equivalent to the diagonal matrix \[D = \mathrm{diag}(2,2,-10,-140,1,-14)\]
    which is $\mathbb{Q}$-equivalent to the standard quadratic form $\mathrm{diag}(1,1,1,-1,-1,-1)$ as follows:

    \begin{align*}
        \mathscr{A} &\simeq_{\mathbb{Q}}
         \begin{bmatrix*}[r]
        2 & 0 & 0 & -4 & 0 & 0 \\
        0 & 2 & -3 & 0 & 0 & 0 \\
        0 & -3 & 2 & -1 & 0 & 0 \\
        -4 & 0 & -1 & 2 & 0 & 0 \\
        0 & 0 & 0 & 0 & 9/4 & -135/4 \\
        0 & 0 & 0 & 0 & -135/4 & 9/4
\end{bmatrix*}\\
&\simeq_{\mathbb{Q}} \mathrm{diag}(2,2,-10,-140,1,-14)\\
&\simeq_{\mathbb{Q}} \mathrm{diag}(1,1,1,-1,-1,-1)
\end{align*}
which can be verified by Hasse--Minkwoski invariants.
Thus there exists a matrix $P$ in $\mathrm{GL}(6,\mathbb{Q})$ such that $\mathscr{A}=PDP^T$.
The diagonal quadratic matrix $D$ is preserved by
    \[c(t)=
\left[
\begin{array}{c|cc}
  I_{4} & \multicolumn{2}{c}{} \\
  \hline
   & \cosh t & \frac{1}{\sqrt{14}}\sinh t \\
    & \frac{1}{\sqrt{14}}\sinh t & \cosh t
\end{array}
\right]
    \]    
for every real numbers $t$. Let $\rho$ be the representation defined by $\rho(\gamma)=P^{-1}\rho_{\mathscr{A}}(\gamma)P$ for every $\gamma\in W$, where $\rho_{\mathscr{A}}$ is the representation introduced in Example~\ref{geom_repn}. 
\begin{figure}[!t]
\centering
\begin{tikzpicture}[
  line width=0.7pt,
  edge label/.style={font=\normalsize,inner sep=2pt},
  face/.style={font=\large},
  sigma boundary/.style={blue!65!black,line width=1pt}
]
  \coordinate (A) at (-2, 2);
  \coordinate (B) at ( 2, 2);
  \coordinate (C) at ( 2,-2);
  \coordinate (D) at (-2,-2);

  \coordinate (a) at (-0.65, 1.15);
  \coordinate (b) at ( 0.65, 1.15);
  \coordinate (c) at ( 0.65,-0.15);
  \coordinate (d) at (-0.65,-0.15);

  \coordinate (T12) at (-0.65,0.5);
  \coordinate (T13) at ( 0.65,0.5);
  \coordinate (T34) at ( 2.00,0.0);
  \coordinate (T24) at (-2.00,0.0);

  \fill[blue!12]
    (T12)--(T13)--(T34)--(T24)--cycle;

  \draw (A) -- node[edge label,above] {$2$} (B);
  \draw (B) -- node[edge label,right,pos=0.3] {$3$} (C);
  \draw (C) -- node[edge label,below] {$2$} (D);
  \draw (D) -- node[edge label,left,pos=0.7] {$2$} (A);

  \draw (a) -- node[edge label,above] {$3$} (b);
  \draw (b) -- node[edge label,right,pos=0.75] {$2$} (c);
  \draw (c) -- node[edge label,below] {$3$} (d);
  \draw (d) -- node[edge label,left,pos=0.25] {$2$} (a);

  \draw (A) -- node[edge label,above right] {$3$} (a);
  \draw (B) -- node[edge label,above left]  {$2$} (b);
  \draw (C) -- node[edge label,below left]  {$2$} (c);
  \draw (D) -- node[edge label,below right] {$3$} (d);

  \draw[sigma boundary]
    (T24)--(T12)--(T13)--(T34);

  \draw[sigma boundary,dashed]
    (T34)--(T24);

  \node[blue!65!black,font=\large]
    at (0,0.23) {$\Sigma$};

  \node[face] at ( 0.00, 0.85) {$F_1$};
  \node[face] at (-1.32, 0.90) {$F_2$};
  \node[face] at ( 1.32, 0.90) {$F_3$};
  \node[face] at ( 2.75, 0.00) {$F_4$};
  \node[face] at ( 0.00, 1.65) {$F_5$};
  \node[face] at ( 0.00,-1.05) {$F_6$};

\end{tikzpicture}

\caption{A Coxeter cube with a totally geodesic $2$-orbifold $\Sigma$ bounded by the facets $F_1,F_2,F_3,F_4$. The lower region is labeled \(\mathcal{O}_1\), and the upper region is labeled \(\mathcal{O}_2\).}
\label{fig:cube-sigma}

\end{figure}
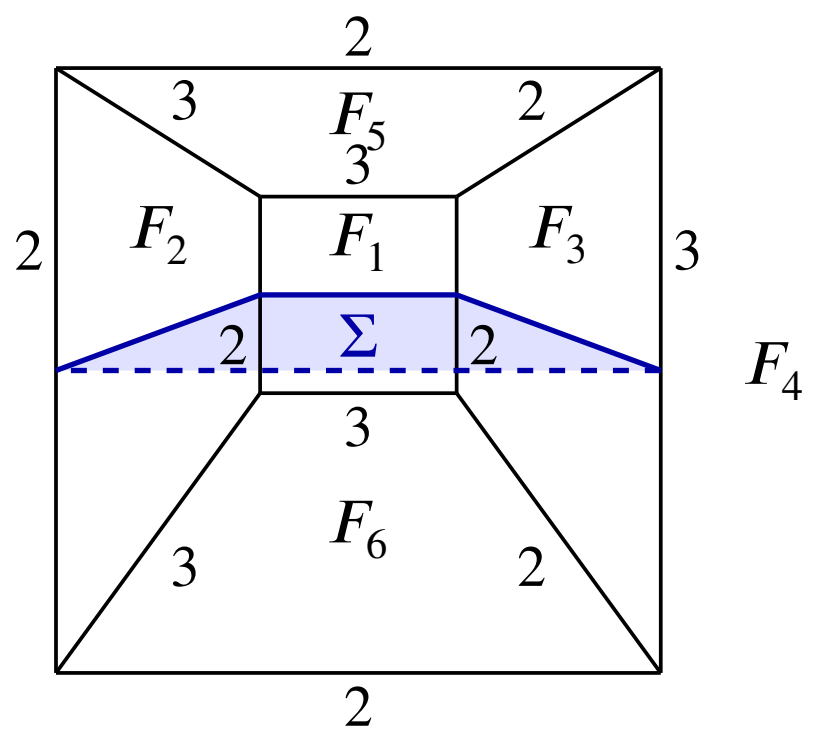
Let $\Sigma$ be a totally geodesic 2-orbifold surrounded by facets $F_1,F_2,F_3$, and $F_4$ (See Figure \ref{fig:cube-sigma}). Since $W_{23}$ is isomorphic to $\pi_1(\mathcal{O}_{23})$, $W_{23}$ can be identified with the amalgamated product $\pi_1^{\text{orb}}(\mathcal{O}_1)*_{\pi_1^{\text{orb}}(\Sigma)}\pi_1^{\text{orb}}(\mathcal{O}_2)$. Put $\gamma_i$ be a reflection in $W_S$ along the face $F_i$. It is routine to check that $\rho(\pi^{\text{orb}}_1(\Sigma))=\langle\rho(\gamma_1),\rho(\gamma_2),\rho(\gamma_3),\rho(\gamma_4)\rangle$ commutes with $c(t)$. Using the matrix $c(t)$, we obtain bending deformations $\rho_t$ of $\rho$:
\begin{align*}
    \rho_t(\gamma)=
    \begin{dcases*}
        \rho_{\mathscr{A}}(\gamma) & for $\gamma\in \pi_1^{\mathrm{orb}}(\mathcal{O}_1)$\\
        c(t)\rho_{\mathscr{A}}(\gamma)c(t)^{-1} & for $\gamma\in \pi_1^{\mathrm{orb}}(\mathcal{O}_2)$
    \end{dcases*}.
\end{align*}
Since $\rho(t)$ varies continuously with $t$, Fact $\ref{thm:HpqCC}$ implies that \[\rho_t\in i_*(\mathrm{Hom}^{\text{cc}}_0(W_S,\mathrm{O}_f(3,3)))\cap \mathrm{Hom}^{\text{ref}}(W_S,\mathrm{GL}(6,\mathbb{R}))\]
for every real number $t$. Thus $\rho_t$ are discrete and faithful representation up to conjugate.
In particular, by solving Pell's equation $x^2-14y^2=1$, we obtain the integral matrix
    
    \[c=
\left[
\begin{array}{c|cc}
  I_{4} & \multicolumn{2}{c}{} \\
  \hline
   & 15 & 56 \\
    & 4 & 15
\end{array}
\right]
    \]
preserving $D$. In other words, $P^{-1}cP$ preserves the symmetric matrix $\mathscr{A}$. Since $P^{-1}cP$ has characteristic polynomial in $\mathbb{Z}[x]$, there is some $k > 0$ so that $(P^{-1}cP)^k \in \mathrm{SL}(6, \mathbb{Z})$; see \cite[Prop 3.2]{long2024zariski},\cite[Lem 4.1]{ZSHORNACK_2025}. Conjugating $\rho_{\mathscr{A}}(\gamma_5)$ with $(P^{-1}cP)^{kn}$ for every integer $n$, we get an infinite sequence $(\rho_n)_n$ in $\mathrm{Hom}(W_{23},\mathrm{O}_f(3,3))$. Hence $(\rho_n(W_{23}))_n$ are the desired family of representations.
\end{proof}

\begin{remark}
\
\begin{enumerate}
\item Using SageMath, we used a brute-force search over all triples $(u,v,w)$ with $u,v,w\leq 50$ and found the triple $(5,3,4)$ for which $\mathscr{A}$ is $\mathbb{Q}$-equivalent to the standard quadratic form. The same approach also gives examples for the other cubes considered in \cite{choi2012projective}.
\item It follows that the family of representations $(\rho_n(W_{23}))_n$ is pairwise nonconjugate, since 
$\mathrm{tr}(\rho_n(\gamma_5)\rho_n(\gamma_6))$ is strictly increasing.
\end{enumerate}
\end{remark}
\begin{lemma}\label{index lemma}
    Let $(\Gamma_n)$ be an infinite sequence of subgroups in a group $G$ each of which is isomorphic to a fixed finitely generated group $\Gamma$. Suppose that $G'$ is commensurable to $G$. Then there exists an infinite sequence $(\Gamma'_k)$ in $G'$ such that $\Gamma'_k$ is isomorphic to a fixed subgroup $\Gamma'$ of finite-index in $\Gamma$.
\end{lemma}
\begin{proof}
    Since $G'$ is commensurable to $G$, we see that
    \[[\Gamma_n:\Gamma_n\cap G']\leq [G:G\cap G']<\infty\]
    On the other hand, $\Gamma$ is finitely generated and thus has only finitely many subgroups of index at most $[G:G\cap G']$. Hence there is an infinite subsequence $(\Gamma_{n_k}\cap G')$ in $G'$ each of which is isomorphic to a fixed subgroup of finite-index in $\Gamma$.
\end{proof}

We can now prove the theorem by Selberg's lemma.

\begin{proof}[Proof of Theorem 1.1]
Let $\Lambda$ be a torsion free subgroup of $\mathrm{SL}(4,\mathbb{Z})$ induced by Selberg's lemma. By Lemma \ref{index lemma}, we see that there is an infinite sequence $(\Gamma'_n)$ in $\varphi(\Lambda)$ such that $\Gamma'_n$ is isomorphic to a fixed subgroup $\Gamma'$ with finite index in $\Gamma$. 

Since $\varphi|_{\Lambda}$ is injective, 
\[\varphi^{-1}(\Gamma'_n)\cap \Lambda \cong \varphi|_{\Lambda}(\varphi^{-1}(\Gamma'_n)\cap\Lambda)=\Gamma'_n\cap \varphi(\Lambda)=\Gamma'_n.\]
Hence there is a closed hyperbolic 3-manifold $M$ such that $\pi_1(M)$ admits infinitely many discrete faithful representations into $\mathrm{SL}(4,\mathbb{Z})$ by applying Selberg's lemma again.

The resulting representation $\rho_t$ is also a representation of $\pi_1^{\text{orb}}(\mathcal{O}_{23})$ as a reflection group by Theorem $\ref{thm:folklore}$, and preserves a nonzero symmetric bilinear form of signature $(3,3)$. Thus it follows from Fact $\ref{fact:zariski_dense}$ that the image of $\rho_t$ is Zariski-dense in $\mathrm{O}(3,3)$.
\end{proof}

\section{Noncommensurablity}
In this section, we shall follow the approach of \cite{10.1215/ijm/1498032023} to prove Theorem 1.2. Although this theorem is a slight modification, we include the whole proof for completeness.

\begin{lemma}
    Suppose that $H$ is a Zariski-dense (closed) Lie subgroup in a simple Lie group $G$. Then $H$ is either discrete or topologically dense in the identity component $G^0$ of $G$. 
\end{lemma}
\begin{proof}
    We may assume that $H$ is topologically closed in $G$. Let $\mathfrak{g}$ and $\mathfrak{h}$ be Lie algebras of $G$ and $H$. Recall that the adjoint map $\mathrm{Ad}:G\rightarrow \mathrm{Aut}(\mathfrak{g})\subseteq \mathrm{GL}(\mathfrak{g})$ is algebraic. Since $\mathrm{Ad}(H)\mathfrak{h}\subseteq \mathfrak{h}$ and $H$ is Zariski-dense, we see that $\mathrm{Ad}(G)\mathfrak{h}\subseteq \mathfrak{h}$, and thus $[\mathfrak{g},\mathfrak{h}]\subseteq\mathfrak{h}$. Simplicity of $G$ gives either $\mathfrak{h}=0$ or $\mathfrak{h}=\mathfrak{g}$. This concludes the result.  
\end{proof}
Given a subset $\Omega$ in $\mathbb{H}^{p,q}$, we denote its stabilizer as
\[\mathrm{Stab}(\Omega):=\{g\in \mathrm{SO}(p,q+1)|\,g\,\Omega=\Omega\}.\]
\begin{lemma}
    Suppose that there is a Zariski-dense group $\Delta$ in $\mathrm{SO}(p,q+1)$ such that it acts on a nonempty proper subset $\mathcal{C}$ of $\mathbb{H}^{p,q}$. Then $\mathrm{Stab}(\mathcal{C})$ is discrete in $\mathrm{SO}(p,q+1)$.
\end{lemma}
\begin{proof}
    Since $\Delta$ is Zariski-dense, the topological closure $\overline{\mathrm{Stab}(\mathcal{C})}$  is also Zariski-dense in $\mathrm{SO}(p,q+1)$. By the previous lemma, $\overline{\mathrm{Stab}(\mathcal{C})}$ is either discrete or topologically dense in $\mathrm{SO}_0(p,q+1)$. If it is in the latter case, $\mathrm{SO}_0(p,q+1)$ must preserves the proper set subset $\mathcal{C}$ of $\mathbb{H}^{p,q}$ by transitivity of $\mathrm{SO}_0(p,q+1)$. This is a contradiction.
\end{proof}

Unlike representations in the $\mathrm{SL}(3,\mathbb{R})$-Hitchin setting treated in \cite{10.1215/ijm/1498032023}, $\mathbb{H}^{p,q}$-convex cocompact representations exhibit distinct behaviors. Since the stabilizer of a point in \(\mathbb{H}^{p,q}\) under \(\mathrm{PO}(p,q+1)\) is isomorphic to the noncompact group \(\mathrm{O}(p,q)\), a discrete subgroup need not act properly discontinuously on \(\mathbb{H}^{p,q}\).
Although a discrete subgroup preserving a properly convex subset \(\mathcal{C}\) with nonempty interior acts properly discontinuously on its interior, this need not extend to all of \(\mathcal{C}\). Moreover, while the convex domain of a Benoist representation, which is in fact a $\mathrm{SL}(3,\mathbb{R})$-Hitchin representation, is uniquely determined by its limit set, an invariant convex domain of an $\mathbb{H}^{p,q}$-convex cocompact representation need not be. We therefore proceed this proof instead with a \(\Gamma\)-invariant complete maximal spacelike \(p\)-submanifold, equipped with its induced Riemannian metric.

\begin{lemma}\label{prop_disc}
    Let $M$ be the complete maximal $p$-submanifold for $\Gamma$ in $\mathbb{H}^{p,q}$. If $\mathrm{Stab}(M)$ is discrete, then $\mathrm{Stab}(M)$ acts properly discontinuously on $M$.
\end{lemma}
\begin{proof}
This is immediate from Corollary 2.5 and Theorem 3.8 in \cite{viaggi2025geometric}.
\end{proof}
\begin{proof}[Proof of Theorem \ref{1.2}]
    Suppose, for the sake of contradiction, that there are infinitely many commensurable subsequences. Its limit sets are $\mathrm{O}(p,q+1)$-equivalent. By conjugating a suitable element, we may assume that its limit sets are the same as the limit set of complete maximal $p$-submanifold $M$.

    Since $\mathrm{Stab}(M)$ contains a $\mathbb{H}^{p,q}$-convex cocompact representation, by Lemma \ref{prop_disc} and Fact \ref{complete_maximal_submanifold}, $M/\mathrm{Stab}(M)$ is isomorphic to the fundamental group of a compact hyperbolic orbifold. Note that the fundamental group of a compact orbifold is finitely generated. Thus it has only finitely many subgroups of a fixed index.
    
    The subgroups $\rho_n(\Delta)$ in $\mathrm{Stab}(M)$ have the finite index $[\mathrm{Stab}(M):\rho_n(\Delta)]$ which equals 
    \[\frac{\mathrm{Vol}(M/\rho_n(\Delta))}{\mathrm{Vol}(M/\mathrm{Stab}(M))}=\frac{||M/\rho_n(\Delta))||}{||M/\mathrm{Stab}(M))||}\]
    by proportionality principle of simplicial volume. The work of Mineyev implies that the simplical volume of a hyperbolic group is positive and independent of $\rho_n$ (see \cite{mineyev01}). It follows that the ratio in the left-hand side is a positive constant.
    Therefore, $\mathrm{Stab}(M)$ contains at least two conjugate subgroups isomorphic to $\rho_n(\Delta)$. Since groups commensurable up to conjugacy have $\mathrm{O}(p,q+1)$-equivalent limit sets, this contradicts our assumption.
\end{proof}
\begin{remark}
    A similar arguement is also applicable to Benoist representations by replacing complete maximal $p$-submanifolds with affine spheres; see \cite{CY76,CY86,viaggi2025geometric}. 
\end{remark}

\newpage
\section{Appendix}
In this appendix we exhibit the concrete computation in Lemma 3.3. 
\subsection{Concrete representations}
Note that the procedure of bending deformation only changes the reflection $\rho(\gamma_5)$ of $F_5$.\\
\begin{minipage}{0.45\textwidth}
\begin{align*}
P&=\begin{bmatrix*}[r]
1 & 0 & 0 & 0 & 0 & 0 \\
0 & 1 & 0 & 0 & 0 & 0 \\
0 & -3/2 & 1/2 & 0 & 0 & 0 \\
-2 & 0 & 1/5 & 1/5 & 0 & 0 \\
-1/2 & -1/2 & 3/10 & 1/20 & 3/2 & 0 \\
-1/2 & -1/2 & 3/10 & 1/20 & -45/2 & 6
\end{bmatrix*}\\
\rho(\gamma_1)&=\begin{bmatrix*}[r]
-1 & 0 & 0 & 0 & 0 & 0 \\
0 & 1 & 0 & 0 & 0 & 0 \\
0 & 0 & 1 & 0 & 0 & 0 \\
4 & 0 & 0 & 1 & 0 & 0 \\
1 & 0 & 0 & 0 & 1 & 0 \\
1 & 0 & 0 & 0 & 0 & 1
\end{bmatrix*}\\
\rho(\gamma_3)&=\begin{bmatrix*}[r]
1 & 0 & 0 & 0 & 0 & 0 \\
0 & 1 & 3 & 0 & 0 & 0 \\
0 & 0 & -1 & 0 & 0 & 0 \\
0 & 0 & 1 & 1 & 0 & 0 \\
0 & 0 & 0 & 0 & 1 & 0 \\
0 & 0 & 0 & 0 & 0 & 1 
\end{bmatrix*}\\
\rho(\gamma_5)&=
\begin{bmatrix*}[r]
1 & 0 & 0 & 0 & 1 & 0 \\
0 & 1 & 0 & 0 & 1 & 0 \\
0 & 0 & 1 & 0 & 0 & 0 \\
0 & 0 & 0 & 1 & 0 & 0 \\
0 & 0 & 0 & 0 & -1 & 0 \\
0 & 0 & 0 & 0 & 34 & 1
\end{bmatrix*}\\
\rho_n(\gamma_5)&=(PcP^{-1})^{-2n}\rho(\gamma_5)(PcP^{-1})^{2n}
\end{align*}
\end{minipage}
\hfill
\begin{minipage}{0.45\textwidth}
\begin{align*}
c &= \left[\begin{array}{rrrr|rr}
1 & 0 & 0 & 0 & 0 & 0 \\
0 & 1 & 0 & 0 & 0 & 0 \\
0 & 0 & 1 & 0 & 0 & 0 \\
0 & 0 & 0 & 1 & 0 & 0 \\
\hline
0 & 0 & 0 & 0 & 15 & 56 \\
0 & 0 & 0 & 0 & 4 & 15
\end{array}\right]\\
\rho(\gamma_2)&=\begin{bmatrix*}[r]
1 & 0 & 0 & 0 & 0 & 0 \\
0 & -1 & 0 & 0 & 0 & 0 \\
0 & 3 & 1 & 0 & 0 & 0 \\
0 & 0 & 0 & 1 & 0 & 0 \\
0 & 1 & 0 & 0 & 1 & 0 \\
0 & 1 & 0 & 0 & 0 & 1
\end{bmatrix*}\\
\rho(\gamma_4)&=\begin{bmatrix*}[r]
1 & 0 & 0 & 4 & 0 & 0 \\
0 & 1 & 0 & 0 & 0 & 0 \\
0 & 0 & 1 & 1 & 0 & 0 \\
0 & 0 & 0 & -1 & 0 & 0 \\
0 & 0 & 0 & 0 & 1 & 0 \\
0 & 0 & 0 & 0 & 0 & 1
\end{bmatrix*}\\
\rho(\gamma_6)&=\begin{bmatrix*}[r]
1 & 0 & 0 & 0 & 0 & 1 \\
0 & 1 & 0 & 0 & 0 & 1 \\
0 & 0 & 1 & 0 & 0 & 0 \\
0 & 0 & 0 & 1 & 0 & 0 \\
0 & 0 & 0 & 0 & 1 & 34 \\
0 & 0 & 0 & 0 & 0 & -1
\end{bmatrix*}\\
\end{align*}
\end{minipage}

\medskip

\bibliographystyle{alpha}
\bibliography{citation}
\end{document}